%%%%%%  THIS IS A PLAIN TEX_FILE.

%      HAUSMANN--KNUTSON
%  POLYGON SPACES AND GRASSMANNIANS ---  July 13 1995

%STYLE
\catcode`\@=11
\magnification=1200
\hsize=140mm        \vsize=200mm
\hoffset=0mm         \voffset=0mm
\pretolerance=500     \tolerance=1000    \brokenpenalty=5000
\parindent=0mm
\footline{\ifnum\pageno=1\hfill \else\hfill -\quad\folio\quad - \hfill\fi}
%  FONTS
%\font\tmX=times at 10pt
\def\tm{}
%%\font\tmXIIbf=timesB at 12pt
%%\font\tmXit=timesI at 10pt
%\font\tmXbf=timesB at 10pt
%\font\tmVII=times at 7pt
%\font\tmV=times at 5pt
%\def\tm{\fam0\tmX}  \tm
%\def\it{\fam\itfam\tmXit}
%\def\bf{\fam\bffam\tmXbf}
%\font\bbX=msym10  at 10pt
%\font\bbVII=msym7 at 7pt
%\font\bbV=msym10 at 5pt
%\newfam\bbfam  \textfont\bbfam=\bbX   \scriptfont\bbfam=\bbVII
%\scriptscriptfont\bbfam=\bbV  \def\bb{\fam\bbfam\bbX}
%\font\mbfX=cmbxti10  at 10pt
%\font\mbfVII=cmbxti10 at 7pt
%\font\mbfV=cmbxti10 at 5pt
%\newfam\mbffam  \textfont\mbffam=\mbfX   \scriptfont\mbffam=\mbfVII
%\scriptscriptfont\mbffam=\mbfV  \def\mbf{\fam\mbffam\mbfX}

\def\calc{{\cal C}}  
  
 \def\calh{{\cal H}} \def\calm{{\cal M}}

%\font\sc=cmcsc10
\def\sc{\sl}
%  MACROS

\def\hfb{\hfil\break}
\def\SS{\S \kern .2em}
 \def\decale#1{\par\noindent\hskip 3em\llap{#1\enspace}\ignorespaces}
\def\apos#1#2{\footnote{\kern 1pt\up{#1}\kern .6pt}{ #2}}

\countdef\parno=20  \parno=0
\countdef\enno=21   \enno=0
\def\noupar#1{\global\advance\parno by 1 \enno=0 {\bf \the\parno .\quad #1%
\bigskip\relax}}
\def\endproc{\endgroup\relax\medskip}
\def\proc#1{\medskip\global\advance\enno by 1 {\bf (\kern 1pt \the\parno \kern1
 pt.
\kern -2pt \the\enno \kern 1 pt) \kern 3pt #1 \kern 3pt } \begingroup\it  }
\def\preu{\sc Proof: \ \tm}
\def\cote{\medskip\global\advance\enno by 1 {\bf (\kern 1pt \the\parno \kern1
 pt.
\kern -2pt \the\enno \kern 1 pt) \kern 6pt  } }
\def\forcote{\global\advance\enno by 1 (\kern 1pt \the\parno \kern1 pt.
\kern 1pt \the\enno \kern 1 pt)}
\def\key{}        %REFERENCES
\def\auteur{}
\def\titre{} \def\journal{} \def\numero{} \def\annee{} \def\pp{}
\def\ref{\begingroup}
\def\endref{\line{\hbox to 1.6truecm{[\key ]\hfill} {\sc \auteur}
\hfill}\vskip -5pt  \quad {\it \titre} \  \journal \  {\bf \numero} \
\annee  \  \pp \relax \vskip 0.4 truecm \endgroup\goodbreak }

\def\frac#1#2{{#1 \over #2}}
\def\pcirc{\kern .7pt {\scriptstyle \circ} \kern 1pt}

\dimendef\dihfl=30  \dihfl = 10 mm    % DIAGRAMMES
\dimendef\divfl=31  \divfl = 5 mm
\def\fl#1{\buildrel{#1}\over{\longrightarrow}}
\def\hfl#1#2{\smash{\mathop{\hbox to \dihfl {\kern 3pt\rightarrowfill\kern
 3pt}}%
\limits^{\scriptstyle#1}_{\scriptstyle#2}}}
\def\hlfl#1#2{\smash{\mathop{\hbox to \dihfl {\kern 3pt\leftarrowfill\kern
 3pt}}%
\limits^{\scriptstyle#1}_{\scriptstyle#2}}}
\def\vfl#1#2{\llap{$\scriptstyle#1$}\left\downarrow\vbox to \divfl{}
\right.\rlap{$\scriptstyle#2$}}

\def\diagram#1{\def\normalbaselines{\baselineskip=0pt\lineskip=10pt
 \lineskiplimit=1pt} \matrix{#1}}
\def\cqfd{\unskip\kern 6pt\penalty 500
\raise -2pt\hbox{\vrule\vbox to10pt{\hrule width 4pt\vfill\hrule}\vrule}\par}
\def\up#1{$^\hbox{#1}$}
\def\tvi{\vrule height 12pt depth 5pt width 0pt}

\catcode`\@=12

\def\picture{\line{\hfill  PICTURE \hfill}}
\def\connect{\ {\scriptstyle\bullet}\!\! \rightarrow \!\!\! -
\!\!{\scriptstyle\bullet} \ }
\def\kah{K\"ahler\ }
\def\fon#1#2{{}^{#1}\!{\cal F}^#2}
\def\pol#1#2{{}^{#1}\!{\cal P}^{#2}}
\def\tpol#1#2{{}^{#1}\!\widetilde{\cal P}^{#2}}
\def\bbr{{{\bf R}}} \def\bbz{{\bf Z}} \def\bbc{{\bf C}}
\def\bbh{{{\bf H}}}
\def\h{\, \breve{}}
\def\apos#1#2{\footnote{\up{#1}\kern .6pt}{\kern 0 truecm \tmVII #2}}
\def\pair#1#2{\langle #1, #2\rangle}
\def\stic{{\bf V}_2(\bbc^m)} \def\stir{{\bf V}_2(\bbr^m)}
\def\grc{{\bf G}_2(\bbc^m)} \def\grr{{\bf G}_2(\bbr^m)}
\def\ogrr{{\bf \tilde G}_2(\bbr^m)}

\parskip= 3truemm plus .2 mm minus .2 mm
\def\preu{{\sc Proof: \ }}
\def\cit#1{(\the\parno.\the\count#1)}
\def\ecit#1#2{(#1.\the\count#2)}
\def\aecit#1#2{(\the\count#1.\the\count#2)}

%%%%%  BEGIN DOCUMENT

\null\vskip 1.3 cm plus .2 cm
\centerline{{\bf
  POLYGON SPACES AND GRASSMANNIANS}}

\vskip 0.5 truecm

\centerline{{\bf Jean-Claude HAUSMANN}
\footnote{$^1$}{University of Geneva, Switzerland; hausmann@ibm.unige.ch}
{\bf and Allen KNUTSON}
\footnote{$^2$}{Massachusetts Institute of Technology;
 allenk@alumni.caltech.edu}
}

\vskip 0.5 truecm
\line{\hfill July 1995 \hfill}

\vskip 1.5 truecm plus .2 cm minus .2 cm
\setbox25=\vtop{\hsize 14 truecm
{\sc Abstract}. We study the moduli spaces of polygons in $\bbr^2$ and $\bbr^3$,
identifying them with subquotients of 2-Grassmannians
using a symplectic version of the Gel$'$fand-MacPherson correspondence.
We show that the bending flows defined by Kapovich-Millson
arise as a reduction of the Gel$'$fand-Cetlin system on the
Grassmannian, and with these determine the pentagon and hexagon spaces
up to equivariant symplectomorphism.
Other than invocation of Delzant's theorem,
our proofs are purely polygon-theoretic in nature.}
\line{\hfill \box25 \hfill }
\vskip 0.7 truecm

\noupar{Introduction}

Let $\tpol{m}{k}$ be the space of $m$-gons in $\bbr^k$ up to translation
and
positive homotheties. This
space comes with several  structures: an action of $O(k)$, an action
of $S_m$
permuting the edges, and a function $\ell : \tpol{m}{k}\fl{}{} \bbr^m$
taking a
polygon $\rho$ to the lengths of its
edges (once the perimeter of $\rho$ is fixed).
The quotients of $\tpol{m}{k}$ by $SO_k$ (or $O_k$) are the moduli
spaces
$\pol{m}{k}_+$ (respectively, $\pol{m}{k}$).
Fixing a reflection in $O(k)$ provides an involution on $\tpol{m}{k}$
and $\pol{m}{k}_+$ whose fixed point sets are $\tpol{m}{{k-1}}$ and
$\pol{m}{k-1}$.
The goal of this paper is to
understand the topology of these various  spaces and the geometric
structures that
they naturally carry when $k=2$ or $3$. They are closely related to more
familiar
objects
(Grassmannians, projective spaces, Hopf bundles,  etc.).
The spaces $\pol{m}{k}(\alpha ):=\ell^{-1}(\alpha )$
of polygons with given side-lengths $\alpha \in \bbr^m$
are of particular interest.

The great miracle occurs when $k=3$, because $\bbr^3$ is isomorphic to the space
$I\bbh$ of pure imaginary quaternions, and the 2-sphere in $\bbr^3$ is
K\"ahler. The tools of symplectic geometry can then be used. Most prominent is a
symplectic version of the Gel$'$fand-MacPherson correspondence identifying the
 spaces
$\pol{m}{3}(\alpha )$ as symplectic quotients of the Grassmannian of 2-planes in
$\bbc^m$.

While this paper illustrates many phenomena in symplectic geometry, the
proofs are entirely polygon-theoretic and involve only classical differential
topology. Nonetheless, many of the examples are
new, interesting in their own right and instructive for both fields.
\goodbreak

Among our results:

\decale{1.} The identification of the polygon space $\pol{m}{3}$ with
$\grc/(U(1)^m)$
intertwines complex conjugation on the
complex Grassmannian (with fixed point set the real Grassmannian)
and spatial reflection on the polygon moduli space (with fixed point set
planar polygons). The fact that 3-dimensional
and planar polygons have the same allowed values of $\ell$ is then an
illustration of a theorem of Duistermaat ([Du]). (As is always true,
and yet always mysterious, it is helpful for studying the real case --
here planar polygons -- to extend to the complex case --
here polygons in $\bbr^3$.)

\decale{2.} Identification of the densely defined ``bending flows" [KM2]
on the polygon spaces with the reduction of the Gel$'$fand-Cetlin system
[GS1] on the Grassmannian.

\decale{3.} In some cases, the bending flows are globally defined,
and by Delzant's reconstruction theory the spaces are equivariantly
symplectomorphic to toric varieties
(for instance when $m\leq 6$, as noted in [KM2]).
We give a precise description of the
moment polytope and so explicitly identify the toric varieties.

Contrary to the usual custom,
 it turned out here to be more natural to take
symplectic quotients by {\sl first} quotienting the original manifold
by the group, and to {\sl then} pick out a symplectic leaf of the resulting
Poisson space -- the intermediate
quotient spaces all have natural polygon-theoretic interpretations.
However, they are never complex; readers wishing a more
geometric-invariant-theoretic
construction of these spaces should look at [KM2].

This paper is structured as follows. Section 2 defines the various
flavors of space of polygons. Sections 3 and 4 relate them to Grassmannians,
and prove some facts about the moment map for the torus action
on the Grassmannian by polygon-theoretic means. In section 4 is also
calculated the exact relation between the \kah structures in this
paper and the ones in [KM2]. Section 5 relates the ``bending flows''
of [KM2] with the Gel$'$fand-Cetlin system on the Grassmannian.
Section 6 uses this to calculate the quadrilateral, pentagon and hexagon spaces.
Section 7 lists some open problems.

The relevance of symplectic geometry in polygon theory is an idea of Sylvain
 Cappell
which arose during a visit of first author at
the Courant Institute in Spring 1993.
This author is also grateful to
Lisa Jeffrey and Mich\`ele Audin for useful conversations.
The two authors started this work at the workshop in symplectic geometry
organized in Cambridge by the Isaac Newton Institute (Fall 1994).
The second author would like to thank Richard Montgomery for teaching
him about dual pairs, and Michael Thaddeus for pointing out the
link to moduli spaces of flat connections; also the University
of Geneva for its hospitality while this paper was being written.
\vskip 1 truecm \goodbreak
\vfill\eject

\parno=1
 \noupar{The polygon spaces}   %%%%SECTION 2 NEW

\cote Let $V$ be a real vector space and $m$ a positive integer. Let
 $\fon{m}{{}}(V)$ be the real vector
 space of
all maps $\rho : \{1,2,\dots m\} \to V$ such that $\sum_{j=1}^m\rho (j)=0$.
 An element $\rho  \in \fon{m}{{}}(V)$ will be regarded as a closed polygonal
 path in
 $V$
$$0 \connect \rho (1) \connect \rho (1) + \rho (2) \connect \cdots \connect
\sum_{j=1}^{m}\rho (j) = 0$$
of $m$ steps,
or, alternately, as a configuration  in $V$  (up to translation) of a
polygon of $m$ sides. We shall call an element $\rho  \in \fon{m}{{}}(V)$
an {\it $m$-polygon (in $V$)}  and a {\it proper polygon} when  $\rho (j) \ne 0
 \
\forall
 j$. We use the notation $\fon{m}{k}$ for the space $\fon{m}{{}}(\bbr^k)$.

The group $\bbr_{+}$ of positive homotheties of $V$ acts freely and
 properly on
$\fon{m}{{}}(V)- \{0\}$. The quotient
$\tpol{m}{{}}(V):=(\fon{m}{{}}(V)- \{0\})/\bbr_{+}$
then inherits a structure of smooth manifold diffeomorphic to the unit sphere in
 $V$.
 For instance, $\tpol{m}{k}:=(\fon{m}{k}-\{0\})/\bbr_{+}$ is diffeomorphic to
 the
sphere $S^{k(m-1)-1}$.

\cote \count60= \enno Suppose now that $V$ is oriented and is a Euclidean space, namely $V$ is
 endowed
with a scalar product. The group $O(V)$ of isometries of $V$ acts  on
 $\fon{k}{m}$ and
$\tpol{m}{{}}(V)$; we define the moduli spaces
$$\pol{m}{{}}(V)_+:=SO(V)\backslash \tpol{m}{{}}(V) \quad \hbox{ and } \quad
\pol{m}{{}}(V):=O(V)\backslash \tpol{m}{{}}(V)$$
of $m$-polygons in $V$, up to similitude
(where $SO(V)$ is the identity component of $O(V)$). Observe that any orientation preserving teisometry
 $h : V
\fl{\sim}{} \bbr^k$ produces identifications
$$\pol{m}{{}}(V)_+ \simeq \pol{m}{k}_+ := SO_k\backslash \tpol{m}{k}
\quad \hbox{ and } \quad
\pol{m}{{}}(V) \simeq \pol{m}{k}:= O_k\backslash \tpol{m}{k}$$
We shall use the fact that these identifications do not depend on the choice of
 $h$
and thus $\pol{m}(V)_+$ and $\pol{m}(V)$, for any Euclidean space $V$, are
canonically identified with $\pol{m}{k}_+$ and $\pol{m}{k}$.

\cote The ``degree of improperness"  of polygons provides a stratification
$$\emptyset = E_1\tpol{m}{{}}(V) \subset E_2\tpol{m}{{}}(V)
\subset \cdots  \subset E_{m-1}\tpol{m}{{}}(V) \subset E_m\tpol{m}{{}}(V)
 \subset =
\tpol{m}{{}}(V)$$ where
$$E_j\tpol{m}{{}}(V) := \{\rho \in \tpol{m}{{}}(V)\mid \sharp\{s\mid \rho
 (s)=0\}\geq
 m-j\}$$
The ``open stratum"
$E_j\tpol{m}{{}}(V) - E_{j-1}\tpol{m}{{}}(V)$ is a smooth submanifold
of $\tpol{m}{{}}(V)$
of dimension $(j-1)k-1$ if $k={\rm dim\,}V$. The top open stratum
$\tpol{m}{{}}(V)-E_{m-1}\tpol{m}{{}}(V)$, open and  dense in $\tpol{m}{{}}(V)$,
is the space of proper polygons.

As this stratification is $O(V)$-invariant, it projects onto stratifications
$\{E_j\pol{m}{k}_+\}$ and $\{E_j\pol{m}{k}\}$ of the moduli spaces (using the
canonical identifications of \cit{60}). We shall see in \SS 3 that the above
stratifications describe the singular loci of  smooth orbifold structures on the
 spaces
$\tpol{m}{{}}(V)$, $\pol{m}{k}_+$ and $\pol{m}{k}$.

\cote \count35=\enno
The map $\rho \mapsto |\rho |:=\sum_{j=1}^m|\rho (j)|$
 which
associates to a polygon $\rho $ its total perimeter  is a norm on
 $\fon{m}{{}}(V)$.
We denote by
$S(\fon{m}{{}}(V))$ the sphere of radius $2$ for this norm. Each class in
 $\tpol{m}{{}}(V)$
has a unique representative in $S(\fon{m}{{}}(V))$ which gives a topological
 embedding
$\imath : \tpol{m}{{}}(V)\fl{}{} \fon{m}{{}}(V)$ whose image is
 $S(\fon{m}{{}}(V))$.
The image by $\imath$ of   $E_{m-1}\tpol{m}{{}}(V)$ is the subset of
$S(\fon{m}{{}}(V))$ where $S(\fon{m}{{}}(V))$ fails to be a smooth submanifold
 of
 $\fon{m}{{}}(V)$.
However,
the  restriction of $\imath$
to  each  $E_j\tpol{m}{{}}(V) - E_{j-1}\tpol{m}{{}}(V)$ is a smooth embedding.

The map $\tilde \ell : \fon{m}{{}}(V) \to \bbr^m $ defined by $\tilde \ell(\rho
 ):=
(|\rho (1)|,\dots ,|\rho (m)|)$ associates to a polygon its side-lengths. We
 define
$\ell : \tpol{m}{{}}(V)\fl{}{} \bbr^m$ by $\ell := \tilde \ell \pcirc \imath$.
We shall also use the notation $\ell_i(\rho )$ for $|\rho (i)|$.
These maps
are invariant under the $O(V)$-action and thus define maps (always called
 $\ell$)
$$\ell : \pol{m}{k}_+ \fl{}{} \bbr^m \ \hbox{ and }
\ell : \pol{m}{k} \fl{}{} \bbr^m$$
which are smooth on each open stratum.

\cote Let $\tau : V \fl{}{} V$ be the orthogonal reflection through some
hyperplane $\Pi$ in $V$. One has the involution
$\rho \mapsto \breve\rho
:=\tau \pcirc \rho $
on $\fon{m}{{}}(V)$ and $\tpol{m}{{}}(V)$ whose fixed-point space is
naturally $\fon{m}{{}}(\Pi)$ and $\tpol{m}{{}}(\Pi)$. If $h \in SO(V)$, one has
$$\tau \pcirc h=\underbrace{(\tau \pcirc h\pcirc \tau \pcirc h^{-1})}_{\in
 SO(V)}
\pcirc h\pcirc \tau $$
Hence the involution descends to an involution (still denoted  $\rho \mapsto
\breve \rho  $)  on $\pol{m}{k}_+$. If $\tau'$ is an orthogonal reflection with
respect to another hyperplane $\Pi'$, then the formula
$\tau \pcirc \rho' = (\tau'\pcirc \tau)\pcirc \tau \pcirc \rho$ shows that the
induced involution on $\pol{m}{k}_+$ does not depend on the choice of $\tau$.
The fixed point space of $\,\breve {}\,$ is $\pol{m}{{k-1}}$. Observe that
 $\breve \rho  =
\rho $ in $\pol{m}{k}$. \smallskip

{\sc Examples: } \smallskip

\cote  \count32=\enno {\it Polygons in the line: } \ The space
$\tpol{m}{1}=\tpol{m}{1}_+$  is diffeomorphic to the sphere  $S^{m-2}$. Under
 this
identification,  the $O_1$-action  becomes the antipodal map and thus
 $\pol{m}{1}$ is a
smooth manifold diffeomorphic to  $\bbr P^{m-2}$.  For example,
 $\tpol{3}{1}\simeq S^1$
and $\pol{3}{1}\simeq \bbr P^1$. The stratum $E_2\tpol{3}{1}$ consists of 3
 pairs of
antipodal points and thus  $E_2\pol{3}{1}$ is a set of 3 points, the three
 triangles
with one side of length $0$. This corresponds to the fact that $S(\fon{3}{1})$
 is a
regular hexagon and  $O_1\backslash S(\fon{3}{1})$ is a triangle. Actually, the
map  $\ell : \pol{3}{1} \fl{}{} \bbr^3$ produces a homeomorphism
$$\pol{3}{1}\hfl{\imath}{\simeq } O_1\backslash S(\fon{3}{1})
\hfl{\ell}{\simeq} \{(x,y,z) \in \bbr_{\geq
0}^3\mid x+y+z=2,  \hbox{ and } \pm x\pm y\pm z = 0\}$$

\cote   \count40=\enno
 {\it Polygons in the plane: } \
Identifying $\bbr^2$ with $\bbc$, the space $\fon{m}{2}$ is a complex vector
 space
isomorphic to $\bbc^{m-1}$ and the (free) $SO_2$-action corresponds to the
 diagonal
$U_1$-action. As in \cit{32} one establishes  the diffeomorphisms
$$\diagram{\tpol{m}{2}  & \hfl{\simeq}{}  & S^{2m-3} \cr
\vfl{}{}  &  & \vfl{}{} \cr
\pol{m}{2}_+  & \hfl{\simeq}{}  & \bbc P ^{m-2} \cr}$$
The above diffeomorphisms
conjugate the involutions $\check{}$ with the complex conjugations of
$\bbc^{m-1}$ and $\bbc P ^{m-2}$.
Also, the involution $\check{}$ on $\pol{m}{2}_+$
coincides with the residual $O_1$ action and
therefore $\pol{m}{{2}}$ is the quotient
of $\bbc P ^{m-2}$ by its complex conjugation.

For example, $\tpol{3}{2}$, the space of planar triangles, is
diffeomorphic to the sphere $S^3$. The singular stratum $E_2(\tpol{3}{1})$ is a
 link of
three circles which are $SO_2$-orbits (therefore, any two of them constitute a
 Hopf
link). The quotient $\pol{3}{2}_+$ is identified with $\bbc P^1$ and
$E_2(\tpol{3}{1}_+)$ is a set of three points in $\bbc P^1$. Finally,
$\pol{3}{2}\simeq \bbc P^1/\{z\sim \bar z\}$ is homeomorphic, via the
 length-side map
$\ell$, to the solid triangle  $$\pol{3}{2}=\pol{3}{3} \hfl{\ell}{\simeq }
\{(x_1,x_2,x_3)\in \bbr^3\mid x_1+x_2+x_3=2 \hbox{ and } 0 \leq x_i \leq 1\}$$
with boundary $\pol{3}{1}$.

\vskip 1 truecm \goodbreak

\noupar{Quaternions, Grassmannians and structures on the full polygon spaces}
 %%% SECTION 3

\cote \count43=\enno Let $\bbh= \bbc \oplus \bbc\, j$ be the skew-field of
 quaternions;
the space $I\bbh$ of pure imaginary quaternions is equipped with the
orthonormal basis $i$, $j$ and $k=ij$,  giving rise to an isometry with $\bbr^3$
which turns the pure imaginary part of the quaternionic multiplication $pq$ into
 the usual cross
product $p\times q$.  The space  $\fon{m}{3}$
is thus identified with $\fon{m}{{}}(I\bbh)$ which gives rise to the canonical
identifications on the the various moduli spaces (see \ecit{2}{60}).

Recall that the  correspondence
$$ \eta : u+vj \mapsto \pmatrix{u & v\cr -\bar v & \bar u\cr}$$
gives  an injective $\bbr$-algebra homorphism $\eta : \bbh \fl{}{}
 \calm_{(2\times
2)}(\bbc)$. This
enables the matrix $P \in U_2$ to act on the right or on the left on $\bbh$. It
 also
identifies the group $S^3=U_1(\bbh)$ of unit quaternions with $SU_2$.

\cote \count70=\enno The Hopf map $\phi :
\bbh\fl{}{} I\bbh$  defined by
$$\phi (q)\, := \, \bar q \, i\,  q$$
sends the 3-sphere of radius $\sqrt{r}$ in $\bbh$
onto the 2-sphere of radius $r$ in $I\bbh$. (The formulae given in the original
 paper
by Hopf [Ho, \SS 5] actually correspond to the map $q\mapsto \bar q k q$.)
The equality $\phi(q)=\phi
 (q')$ occurs if and only if $q'= e^{i\theta}\, q$. The map $\phi$ satisfies the
equivariance relation $\phi (q\cdot P) = P^{-1}\cdot  \phi (q)\cdot  P$.
Writing $q=u+vj$ with $u,v\in\bbc$, one has
$$\phi(u+vj)= (\bar u- j \bar v ) i (u+vj) = i(\bar u +j \bar v)(u+vj) =
i \big[ ( |u|^2-|v|^2) + 2\bar u v j \big]$$

\cote \count71=\enno Observe that if $q = s + tj$ with $s,t \in \bbr$, then
 $\phi (q) =
i\, q^2$. This complex plane $\bbr \oplus \bbr\, j$ is the fixed point set of
 the
involution  $a+bj \mapsto  \bar a + \bar b j$ which will be used later.
Its image under $\phi$ is $\bbr\,i\oplus \bbr\,k$.
\smallskip

\cote {\bf Remark: } \ $I\bbh$, with the Lie bracket $[p,q]=pq-qp = 2{\rm
Im\,}(pq)$, is the Lie algebra for the group $U_1(\bbh) \simeq
SU_2\simeq S^3$. The pairing $(q,q')\mapsto -{\rm Re\,}(qq')=\pair{q}{q'}$
identifies $I\bbh$ with its dual. If $\bbh\simeq  \bbc\oplus \bbc$ is endowed
 with the
standard \kah form, then the map $\frac{1}{2}\phi$ is the moment map for
the Hamiltonian action of $U_1(\bbh)$ on $\bbh$ (the factor $\frac{1}{2}$ can be
checked by restricting the action to the $S^1$-action on $\bbc$).
 \count72=\enno  \smallskip

\cote  \count44=\enno Let  $\stic$ be the space of $(m\times 2)$-matrices
$$(a,b)\ :=\ \pmatrix{a_1 & b_1\cr \vdots & \vdots \cr a_m & b_m\cr}\ \in \
\calm_{m\times 2}(\bbc)$$
such that $|a|=|b|=1$ and $\pair{a}{b} = 0$. $\stic$ is the Stiefel
manifold of orthonormal $2$-frames in $\bbc^m$.  The group $U_m$ acts
 transitively
on the
left on $\stic$ producing the diffeomorphism $\stic = U_m/U_{m-2}$.
One has the conjugation on $\stic$ given by
$(a,b) \mapsto (\bar a,\bar b)$
with fixed-point space the Stiefel manifold
 $\stir= O_m/O_{m-2}$  of orthonormal $2$-frames in $\bbr^m$.
Finally, the embedding  $\stic \subset \bbh^m$ given by $(a,b)\mapsto (\dots ,
a_r+b_r\, j,\dots )$ intertwines the conjugation on $\stic$ with the
involution of \cit{71} on $\bbh^m$. One thus gets an embedding
$\stir\subset (\bbr\oplus \bbr j)^m$.

Using the Hopf map $\phi$ of \cit{70},
one defines  the smooth map $\Phi : \stic
\fl{}{} \fon{m}{{}}(I\bbh)\simeq \fon{m}{3}$  by the formula
$$\Phi(a,b) :=
\ \big( \phi (a_1+b_1 \, j), \phi (a_2+b_2\, j),\dots ,\phi (a_m+b_m \,
 j)\big)$$
The fact that $\sum \phi (a_r+b_r \, j)=0$ is equivalent to $\pair{a}{b} = 0$
 and
$|a|=|b|$. As
$|a|=|b|=1$,  the image of $\Phi $ is exactly $S_2(\fon{m}{3})$.
By composing with the projection $\fon{m}{3}-\{0\}\fl{}{} \tpol{m}{3}$, one gets
 a
surjective smooth map $\Phi : \stic \fl{}{} \tpol{m}{3}$.
One checks that  $\Phi (a,b)=\Phi
(a',b')$ if and only if $(a,b)$ and $(a',b')$ are in the same orbit under the
 action
of the maximal torus $U_1^m$ of diagonal matrices in $U_m$. This action is free
 when
none of the $(a_i,b_i)$'s vanishes, namely if and only if $\Phi (a,b)$ is a
 proper
polygon.
As $\Phi (\bar a,\bar b) = \Phi (a,b)\check{} $, the restriction of
$\Phi $ to the fixed points gives a smooth map $\Phi _\bbr : \stir \fl{}{}
 \tpol{m}{{}}(\bbr\, i\oplus \bbr\, k) \simeq \tpol{m}{2}$
with analogous properties.
We  have thus proved

\proc{Theorem}  a)\ The smooth map $\Phi : \stic
\fl{}{} \tpol{m}{3}$ induces a homeomorphism $\widehat \Phi  :
 U_1^m\backslash\stic
\fl{\simeq}{} \tpol{m}{3}$ such that $\widehat \Phi (\bar a,\bar b) = \Phi
(a,b)\h $. The restriction of $\Phi $ above the space of proper polygons
is a smooth principal $U_1^m$-bundle.
\decale{b)}  The smooth map $\Phi_\bbr : \stir \fl{}{}
\tpol{m}{2}$ induces a homeomorphism \hfb $\widehat \Phi_\bbr  :
 O_1^m\backslash\stir
\fl{\simeq}{} \tpol{m}{2}$.
The restriction of $\Phi_\bbr $ above the space of proper
planar polygons  is a principal $O_1^m$-covering.  \endproc  \count30=\enno

\proc{Corollary}  $\tpol{m}{3} \simeq U_1^m\backslash U_m/U_{m-2}$ and \
$\tpol{m}{2} \simeq O_1^m\backslash O_m/O_{m-2}$.
\endproc  \count31=\enno

\cote  Let $\grc$ be the Grassmann manifold of 2-planes in $\bbc^m$.

The map $\stic \fl{}{} \grc$ which associates to $(a,b)$ the plane generated by
 $a$ and
$b$ is the projection $\stic \fl{}{} \stic / U_2$ (a principal $U_2$ bundle),
 for the
natural right action of $U_2$ on $\stic \subset \calm_{m\times 2}(\bbc)$. This
projection is $U_m$-equivariant, equivalent to the projection $U_m/U_{m-2}
\fl{}{} U_m/U_2\times U_{m-2}$.

The map $\Phi : \stic \fl{}{} \tpol{m}{3}$ satisfies
$$\Phi \big ( (a,b) P \big) = P^{-1}\, \Phi  \big (a,b)\, P  \quad \hbox{ for }
\quad  (a,b) \in \stic ,\ P \in U_2\, $$
The conjugation by $P$ being an element of $SO(I\bbh)$, one thus gets
a map (still called $\Phi $) from $\grc$ onto $\pol{m}{3}_+$.
The space $\pol{m}{3}_+$ has a smooth structure on the open-dense subset of
 non-lined
polygons (which is where the $SO_3$-action was free)
 and,
above this open-dense subset, the new map $\Phi$ is smooth.
The map $\Phi$ intertwines  the  involutions and so restricts to  a map
$\Phi_\bbr :  \grr \fl{}{} \pol{m}{2}$, where $\grr$ is the Grassmannian of
2-planes in $\bbr^m$. In this case,  an intermediate object is the Grassmannian
$\ogrr= SO_m/SO_2\times SO_{m-2}$ of
oriented 2-planes in $\bbr^m$ with the smooth map $\Phi_\bbr \ogrr \fl{}{}
\pol{m}{2}_+\simeq \bbc P^{m-2}$. The action of $U_1^m$ on $\stic$ descends to
 an
action on $\grc$ which is no longer effective: its kernel is the diagonal
 subgroup
$\Delta $ of $U_1^m$, the center of $U_m$, isomorphic to $U_1$. The same
holds true in the real case, replacing $U_1$ by $O_1$ (the diagonal subgroup of
$O_1^m$ is also denoted by $\Delta $).

Using Theorem \cit{30}, the reader will easily prove the folowing

\proc{Theorem}  a)\ The  map $\Phi : \grc
\fl{}{} \pol{m}{3}$ induces a homeomorphism \hfb $\widehat \Phi  :
U_1^m\backslash\grc \fl{\simeq}{} \pol{m}{3}$ such that $\widehat \Phi (\bar
a,\bar b) = \Phi (a,b)\h $. The restriction of $\widehat \Phi $ above the space
 of
proper non-lined polygons is a smooth principal $(U_1^m/\Delta )$-bundle.
 \decale{b)}  The smooth map
$\Phi_\bbr : \ogrr \fl{}{} \pol{m}{2}_+$ induces a homeomorphism \hfb $\widehat
\Phi_\bbr  : O_1^m\backslash\ogrr \fl{\simeq}{} \pol{m}{2}_+$. It is a smooth
branched covering and,  restricted above the space of proper
polygons,  a principal $(O_1^m/\Delta )$-covering.
\decale{c)}  The map $\Phi_\bbr : \grr \fl{}{} \pol{m}{2}$ induces a
 homeomorphism
\hfb $\widehat \Phi_\bbr  :
O_1^m\backslash\grr \fl{\simeq}{} \pol{m}{2}$. The restriction of $\widehat \Phi
 $ above the space of
proper non-lined polygons is a  principal $(O_1^m/\Delta )$-covering.
\endproc  \count33=\enno

\proc{Corollary}  One has homeomorphisms between the polygon spaces and the
 double
cosets
\decale{a)}
$\pol{m}{3} \simeq U_1^m\backslash U_m/(U_2\times U_{m-2})$
\decale{b)}
$\pol{m}{2}_+ \simeq S(O_1^{m})\backslash SO_m/(SO_2\times SO_{m-2})$.
\decale{c)} $\pol{m}{2} \simeq O_1^{m}\backslash O_m/(O_2\times O_{m-2})$.
\endproc  \count34=\enno

\cote {\it Example: }\ As in \ecit{2}{40} the example of  planar triangles
 ($m=3$
and $k=2$) is interesting. The Stiefel manifold
 ${\bf V}_2(\bbr^3)$ is diffeomorphic to the unit tangent bundle to $S^2$,
in turn
diffeomorphic to $SO_3$. The oriented Grassmannian ${{\bf \tilde G}_2}(\bbr^3)$
 can
be identified with $S^2$ by associating to an oriented plane its unit normal
 vector.
The smooth map $$\Phi_\bbr : S^2\simeq {{\bf \tilde G}_2}(\bbr^3))\fl{}{}
\pol{3}{2}_+ \simeq
S^2$$ is of degree $4$, branched over the $3$ points. This map can be visualized
 as
follows: tesselate $\bbr^2$ with equilateral triangles. Divide
 $\bbr^2$ by
the subgroup of isometries which preserve the tesselation
 and the orientation (it
thus preserves a checkerboard coloring of the triangle tesselation). This
 quotient
is a well known orbifold structure on $S^2$ with three branched points. The
projection $\bbr^2\fl{}{} S^2$ factors through an octahedron with a chess-board
coloring of its faces. The residual map from this octahedron to $S^2$ is our map
$\Phi_\bbr$.

Take the pullback by $\Phi_\bbr$ of the Hopf bundle $S^3\fl{}{} S^2$. One gets a
map of degree $4$ from some lens space $L$ onto $S^3$, which branched locus the
 link
formed by three $SO_2$-orbits. The lens space will be doubly covered by $SO_3$.
 We
thus gets the map
$$\widetilde \Phi : SO_3\simeq {\bf V}_2(\bbr^3)\fl{}{} \tpol{3}{2}\simeq
S^3$$
of degree $8$. Finally, one has ${\bf  G}_2(\bbr^3)\simeq \bbr P^2$ and
 $\Phi_\bbr$ is
the quotient  of $\bbr P^2$ by the action of  $O_1^3$ on each homogeneous
 coordinate.
This quotient is a 2-simplex and one sees again that
 $\pol{3}{2}$
is a solid triangle.

\cote  {\it Orbifold structures: } \ The maps $\widehat\Phi_\bbr$ and
 $\Phi_\bbr$
provide, for  the spaces $\tpol{m}{2}\simeq S^{2m-3}$ and $\pol{m}{2}_+\simeq
 \bbc
P^{m-2}$, a smooth   orbifold structure. Each point has a neighbourhood
 homeomorphic
to an open set  of the quotient of $(\bbr^2)^{s}$ by a subgroup of $O_1^s$,
 where $O_1$
acts on each  $\bbr^2$ via the antipodal map. Observe that the map $\Phi_\bbr$
is a ``small cover" in the sense of [DJ].
The branched loci are
$E_{m-1}\tpol{m}{2}$ and  $E_{m-1}\pol{m}{2}_+$ respectively. As for
 $\pol{m}{2}$ we
have to add the branched locus $\pol{m}{1}$. The generic  points of $\pol{m}{1}$
have a neighbourhood modelled on the quotient of $\bbc^{m-2}$
 by  complex conjugation.

Analogously, the map $\Phi : \grc\fl{}{} \pol{m}{3}$  gives rise, for the space
$\tpol{m}{3}$, to  a smooth {\it complex orbifold structure}. By that we mean
a space locally modelled on the quotient of
$\bbc^{s}$ by a subgroup of $U_1^s$.
We define the space $\calc^\infty(\pol{m}{3})$ of {\it smooth} maps from
$\pol{m}{3}$ to the reals as the subspace of $\calc^\infty(\grc )$ which is
invariant by the action of $U_1^m$.

\cote \count41=\enno {\it Riemannian and Poisson structures: }\ Let $\calh (m)$
 be
the space of Hermitian $(m\times m)$-matrices, identified with  ${\bf u} _m^*$
 via the
pairing
$$\calh (m)\times {\bf u _m}\fl{}{} {\bbr} \qquad (H,X)\mapsto \frac{i}{2} {\rm
tr}(HX)$$
This identification turns the co-adjoint action of $U_m$ into the  conjugation
action on $\calh (m)$.
Consider the map $\widetilde\Psi : \calm_{m\times 2}(\bbc)\fl{}{} \calh (m)$
 given by
$\widetilde\Psi(a,b):=(a,b)\cdot (a,b)^*$. One has  $\widetilde\Psi(Q\cdot
 (a,b)\cdot P)=
Q\cdot \widetilde\Psi((a,b)\cdot Q^*$ for $P \in U_2$ and $Q\in U_m$ and thus
$\calc:=\widetilde\Psi(\stic )$ is the $U_m$-orbit through ${\rm
 diag}(1,1,0,\dots ,0)$.
This proves that $\widetilde\Psi$ descends to a diffeomorphism $\Psi :
 \grc\fl{\simeq
}{}  \calc$.

The complex vector space $\calm_{m\times 2}(\bbc)$ is endowed with its
classical Hermitian structure $\pair{A}{B}:={\rm tr}(AB^*)$, with associated
 symplectic
form $\omega (\, ,\, )=-{\rm Im} \pair{\, }{\, }$.
The map $\widetilde\Psi$ above and the map $\widetilde \Phi : \calm_{m\times
 2}(\bbc)
\fl{}{} \calh_0 (2)$ given by
$$\widetilde \Phi (a,b) := (a,b)^*\cdot (a,b) - \pmatrix{1 & 0\cr 0 & 1\cr}$$
 are moment
maps for the Hamiltonian actions of $U_m$ and $U_2$ respectively. One has
 $\stic=
\widetilde \Phi ^{-1}(0)$ and thus $\grc$ occurs as symplectic reduction of the
hermitian vector  space $\calm_{m\times 2}(\bbc)$ and thereby inherits a
 $U_m$-invariant
\kah structure, using, for instance [Ki], \SS 1.7. (Strictly speaking, one deals
 in
[Ki] with compact \kah manifolds;
to fulfill this condition, one can first divide $\calm_{m\times 2}(\bbc)-\{0\}$
 by the
diagonal action of $\bbc^*$ to put onself into a complex projective space).
The residual map
$\Psi :  \grc\fl{\simeq }{}  \calc \subset \calh (m) $ is a moment map for the
 action of
$U_m$ on $\grc$.

Being thus a \kah manifold, $\grc$ is a Riemanniann Poisson manifold. This
 structure
descends to the complex orbifold $\pol{m}{3}$: the algebra
$\calc^\infty(\pol{m}{3})$
admits a unique Lie bracket so that the projection
$\grc\fl{}{} \pol{m}{3}$ is a Poisson map.

\vskip 1 truecm \goodbreak

\noupar{Polygons with given sides -- \kah structures}  %%% SECTION 4

We now use the map  $\ell : \tpol{m}{k},\ \pol{m}{k}_+,\ \pol{m}{k} \to \bbr^m $
 defined in
\ecit{1}{35}. Recall that $\ell (\rho )$, for $\rho \in\tpol{m}{k}$, is the
 length of
the successive sides of a representative of $r$ with total perimeter $2$.

For $\alpha =(\alpha _1,\dots ,\alpha _m)\in \bbr_{\geq 0}^m$ with $\sum_{i=1}^m
\alpha _i=2$, we define $$\tpol{m}{k}(\alpha ):=:\tpol{{}}{k} (\alpha ):=\{\rho
 \in
\tpol{m}{k}\mid \ell(\rho )=\alpha )\}\subset \tpol{m}{k}$$
The space $\tpol{{}}{k}(\alpha )$ is invariant under the action of $O_k$. We
 define
the moduli spaces
$$\pol{{}}{k}_+(\alpha ):=SO_k\backslash \tpol{{}}{k}(\alpha )=
\ell^{-1}(\alpha ) \subset \pol{m}{k}_+$$
and
$$\pol{{}}{k}(\alpha ):=O_k\backslash \tpol{{}}{k}(\alpha )
= \ell^{-1}(\alpha ) \subset \pol{m}{k}$$
The space $\tpol{{}}{1} (\alpha )$ consists of a finite number of points and is
generically empty. We call $\alpha $ {\it generic} if $\tpol{{}}{1} (\alpha )=
\emptyset$.

\proc{Theorem} The map $\mu  := \ell\pcirc \widehat \Phi : \grc \fl{}{} \bbr^m$
 is
a moment map for the action of \ $U_1^m$ on $\grc$. \endproc  \count36=\enno
\count42=\enno
\preu
As seen in \ecit{3}{41}, the moment map $\Psi :\grc \fl{}{} \calh (m)$
for the $U_m$-action on $\grc$ is induced
from $\widetilde\Psi : \calm_{m\times 2}(\bbc)\fl{}{} \calh (m)$ given by
$\widetilde\Psi(a,b):=(a,b)\cdot (a,b)^*$. A moment map $\mu $ for the action of
\ $U_1^m$ is obtained by composing $\Psi$ with the projection $\calh (m) \fl{}{}
\bbr^m$ associating to a matrix its diagonal entries. So, if $\Pi \in  \grc$ is
generated by $a$ and $b$ with $(a,b)\in \stic$, one has
$$\mu (\Pi)= (|a_1|^2+|b_1|^2,\dots ,|a_m|^2+|b_m|^2) = \ell\pcirc \widehat \Phi
(a,b)\kern 2 truecm \hbox{\cqfd} $$

A now classic theorem of Atiyah and
Guillemin-Sternberg [Au, \SS III.4.2] asserts that the image of a moment map for
 a torus
 action is a
convex  polytope (the {\it moment polytope}). The restriction of the moment map
 to
the fixed point set of an anti-symplectic involution has the same image [Du].
In our case, one gets these facts directly:

\proc{Corollary} The moment map $\mu : \grc \fl{}{} \bbr^m$
satisfies $\mu (\grc ) = \mu (\grr ) = \Xi_m$, where $\Xi_m$ is the
 hypersimplex
$$\Xi_m := \{(x_1,\dots x_m)\in\bbr^m\mid 0\leq x_i\leq 1 \ \hbox{ and }\
\sum_{i=1}^mx_i=2\}$$
\endproc  \count37=\enno

\preu One has ${\rm Image} (\mu)={\rm Image} (\ell )$ and it is manifest that
 the
${\rm Image}(\ell )\subset \Xi_m$. A proof that
${\rm Image}(\ell )= \Xi_m$ is actually provided in [KM1], Lemma 1, or [Ha]. We
 give
here however another argument, for the pleasure of  constructing  a
continuous section
$\sigma : \Xi_m \fl{}{} \pol{m}{2}$ of $\ell$. If $m=3$,  we have already
mentioned in
\ecit{2}{40} that $\pol{3}{2}$ is homeomorphic to $\Xi_3$ via the map $\ell$.
Let $\alpha \in\Xi_m$. Define  $\beta _i:=\sum_{j=1}^i\alpha _j$ and
$$r(\alpha ) :={\rm min\,}\{i\mid \beta _i \leq 1 \hbox{ and } \beta _{i+1} \geq
 1\}$$
The numbers $\beta _r,\alpha _r,2-\beta _{r+1}$ form a triple of  $\Xi_3$ and
 are
then the lengths of  a unique triangle $\tau (\alpha )\in \pol{3}{2}$ which can
 be
subdivided in the obvious way to define the element
 $\sigma (\alpha )\in \pol{m}{2}(\alpha )$.

\picture

The continuity of $\sigma $ comes from the fact that if
the map $r$ is  discontinuous at some $\alpha $,
the triangle $\tau (\alpha )$ is then
lined. \cqfd \smallskip

{\bf Remarks: \ } 1) Corollary \cit{37} is also a consequence of our stronger
 result
(5.4).

\decale{2)} The word ``hypersimplex" is introduced in [GM]. Observe that $H$
 is
obtained  by taking the convex hull of the middle point of each edge of a
 standard
$(m-1)$-simplex.
\smallskip

We also obtain the critical values of $\mu $ (compare [Ha]):

\proc{Proposition} The set of critical values of $\mu $ on $\grc$ or $\grr$
 consists of
the union of \decale{a)} the boundary of $\Xi_m$. Its pre-image by $\ell$ in
$\pol{m}{k}_+$ \ ($k = 2$ or $3$) is the space of non-proper polygons
$E_{m-1}\pol{m}{k}_+$. \decale{b)} the inner walls $\{(x_1,\dots x_m)\in\Xi_m
 \mid
\exists\ \varepsilon _i=\pm 1 \hbox{ such that } \sum_{i=1}^m\varepsilon
 _ix_i=0\}$
formed by the
non-generic $\alpha $'s. Their pre-image by  $\ell$ in $\pol{m}{k}_+$ contains
 the set
$\pol{m}{1}$ of lined polygons. \endproc  \count38=\enno

\preu The boundary of $\Xi_m$ consists of course of critical values. They are
 the
 image
of the critical points of $\Psi : \stic \fl{}{} \tpol{m}{3}$ (or those of
$\Psi_\bbr$). The inner walls are the image of the critical points of
$\tpol{m}{3} \fl{}{} \pol{m}{3}_+$, where $SO_3$ (or $O_1$) acts non-freely.
 These
points are the lined polygons. \cqfd \smallskip

We have proven most of the main result of this section:  for generic
and proper $\alpha $, the space $\pol{{}}{3} (\alpha )$  is a  \kah
 sub-quotient
of $\grc$.

\proc{Theorem}  For $\alpha \in {\rm int\,} \Xi_m$ generic,
$\pol{{}}{3}_+ (\alpha )$
is a \kah manifold  isomorphic to the \kah
reduction
$U_1^m\backslash \mu ^{-1}(\alpha )$.
The involution $\h $
is antiholomorphic and $\pol{{}}{2} (\alpha )$ can be seen as the real part of
$\pol{{}}{3}_+ (\alpha )$.
\endproc  \count39=\enno

\preu  By \cit{42}, one has $\pol{{}}{3} (\alpha )= \ell^{-1}(\alpha ) =
U_1^m\backslash \mu ^{-1}(\alpha )$ and we have seen in \ecit{3}{33} that
$\widehat \Phi (\bar a,\bar b) = \Phi (a,b)\h $. \cqfd

We shall now compare the \kah structure obtained on
$\pol{{}}{3}_+ (\alpha )$ from the Grassmannian to
that introduced by
Kapovich-Millson ([KM2], \SS 3).
Using the standard cross product $\times $ and scalar product $\pair{.}{.}$ on
$\bbr^3$,  these authors put on the sphere $S^2_r$ of radius $r$ the complex
 structure
$\tilde J$ defined by
$$\tilde J v := \frac{1}{r}\, x\times v \qquad  (v\in T_xS^2_r)$$
and the
\kah metric
$$\tilde h(u,v) :=  \frac{1}{r} \pair{u}{v} -\frac{i}{r^2}\pair{x}{u\times
 v}\qquad (u,v
\in T_xS^2_r)$$ with associated symplectic form $\tilde \omega (u,v) :=
\pair{\frac{x}{r^2}}{u\times v}$. Let $W(\alpha ):=\prod_{i=1}^mS^2_{\alpha
 _i}$.
 The
map $\beta :W_\alpha \fl{}{} \bbr^3$ defined by $\beta (z_1,\dots
,z_m):=\sum_{i=1}^mz_i$ is the moment map for the diagonal action of $SO_3$ on
$W_\alpha $. The space $\pol{{}}{3}_+ (\alpha )$ thus occurs as the symplectic
reduction  $SO_3\backslash \beta ^{-1}(0)$.

\proc{Proposition} The complex structure $J$ and \kah metric $h$ of \cit{39}
 compare
with those $\tilde J$ and $\tilde h$ of Kapovich-Millson in the following way
$$\tilde J = J \quad \hbox{ and } \quad \tilde h(u,v) = 4\,
 {h(u,v)}$$
\endproc  \count45=\enno

\preu Starting from the Hermitian vector space
$\calm = \calm_{m\times 2}(\bbc)$ one sees that $\pol{{}}{3} (\alpha )$ is
 obtained
by two successive symplectic reductions
$$\grc = \widetilde \Phi^{-1}(0)/U_2 \quad \hbox{ and }\quad
\pol{{}}{3} (\alpha )= U_1^m\backslash \mu ^{-1}(\alpha )$$
(We use the notation of \SS 3). One can perform the reductions
in the reverse order. We
first get
$$U_1^m\backslash \widetilde \Psi ^{-1}(\alpha )=\prod_{i=1}^m\bbc
 P_{\alpha
_i}^1$$
where $\bbc P^1_r$ is the quotient of the 3-dimensional sphere
$\{(u,v)\in \bbc^2\mid |u|^2+|v|^2=r\}$ by the diagonal action of $U_1$.
The moment map $\widetilde \Phi : \calm \fl{}{} \calh (2)$ gives a a moment
 map
(still called $\widetilde \Phi$) from the product of projective spaces into
 $\calh_0
(2)$.
We have seen in \ecit{3}{43} and \ecit{3}{44} that one has a commutative diagram
$$\diagram{
\prod_{i=1}^m\bbc P_{\alpha_i}^1  & \hfl{\prod \phi }{\simeq } &
\prod_{i=1}^mS^2_{\alpha_i} \cr
\vfl{\widetilde \Phi}{}  &  & \vfl{}{\beta } \cr
\calh _0 (2)  & \hfl{\psi^{-1}}{\simeq }  & \bbr^3
}$$

To prove Proposition \cit{45}, it is enough to establish that for all $a\in \bbc
P^1_r$, the tangent map
 $T_a\phi
:T_a\bbc P^1_r \fl{}{} T_{\phi (a)}S^2_r$ satisfies
$$T_a\phi(Jv)=\tilde J T_a\phi(v)
\quad  \hbox{ and }\quad \tilde \omega (T_a\phi (v),T_a\phi (J v)) = 4
 \omega
(v,J v)$$
By $U_2$-equivariance, one can restrict ourselves to $a = [\sqrt{r},0]$. The
tangent space $T_a \bbc P^1_r$ is identified with $\{0\}\times \bbc$ and one can
take $v=(0,1)$ and $J v = (0,i)$. One has
$$\phi (a) = (r,0,0) \ ,\ T_a\phi (v) = (0,2\sqrt{r},0) \ , \
T_a\phi (J v) = (0,0, 2\sqrt{r}) = \tilde J T_a\phi (v)$$
and $ \tilde \omega (T_a\phi (v),T_a\phi (J v))=4$ while $\omega
 (v,Jv)=1$.
\cqfd
\medskip
{\bf Remarks: }
\cote  The results of this section show that the spaces $\pol{{}}{3}_+ (\alpha
 )$
for generic $\alpha $ are the symplectic leaves of the Poisson structure
on the regular part of $\pol{m}{3}_+$ given in \ecit{3}{41}. \smallskip
\cote  If one works in the pure quaternions $I\bbh$, the complex structure
 $\tilde
J$ on  $S^2_r$ becomes
$$\tilde J (v) = \frac{q\, v}{|q|} \quad ,\quad
(v\in T_qS^2_r= I\bbh)$$
The sphere $S^2_r$ is a co-adjoint orbit of $U_1(\bbh)$ and the Hermitian form
$\tilde w$ is the Kirillov--Kostant form (see [Gu, Theorem 1.1]). \smallskip
\cote  The isomorphism between the
symplectic reductions of the Grassmannian $\grc$ and the product of $\bbc P
 ^1$'s that
underlies our results \ecit{3}{33}, \cit{39} and the proof of \cit{45} is a
symplectic version of the Gel$'$fand-MacPherson correspondence ([GM] and
 [GGMS]).
The fact that this isomorphism comes  from two reductions of
$\calm$ is the philosophy of ``dual pairs" (see [Mo] and
the references therein).

\vskip 1 truecm \goodbreak

\noupar{The Gel$'$fand-Cetlin action}  %%%%% SECTION 5

On $\fon{m}{k}$ we have so far defined the length functions $\tilde \ell$
measuring the distances between successive vertices. We now introduce
$\tilde d : \fon{m}{k} \to \bbr^m$, $\tilde d(\rho) = (|\rho(1)|,
|\rho(1)+\rho(2)|, \ldots, |\sum_{i=1}^m \rho(i)|)$,
the lengths of the diagonals connecting the vertices to the origin.
(Only $m-3$ of these functions are new, as $\tilde d(\rho)_1 = \tilde
 \ell(\rho)_1$,
$\tilde d(\rho)_{m-1} = \tilde \ell(\rho)_m$,
and $\tilde d(\rho)_m = 0$. Hereafter we write only $\ell_i, d_i$ and the
$\rho$ is to be understood.)

As with $\tilde \ell$, the function $\tilde d$ descends to  continuous but
only generically smooth functions $d$ on
$\tpol{m}{k}$, $\pol{m}{k}_+$ and $\pol{m}{k}$. It is smooth where no $d_i$
 vanishes,
 that is to say the polygon does not return to the origin prematurely. We call
such a polygon $P$ {\it prodigal} and  call $(\ell(P),d(P))$
a  {\it prodigal value}. The set of prodigal
polygons is  open dense in  $\pol{m}{k}_+$ with complement
of codimension $k$.

For $k=3$, there is in [KM1] introduced an action of a torus $T^{m-3}$ on
 prodigal
polygons; the $i$th circle acts by rotating the section of the polygon formed by
 the
first $i$ edges about the $i$th diagonal. (When that diagonal is length zero,
 there is
no well-defined axis about which to rotate, and indeed the action cannot be
 extended
continuously over this subset.) This action plainly preserves
the level sets of the functions $d$, but more is true:

\proc{Theorem [KM1]} On the subspace of prodigal polygons of
 $\pol{{}}{3}_+(\alpha)$,
the  function $d$ is a moment map for these ``bending flows''. \endproc
\count47=\enno

One important consequence of this is that the torus action also preserves
the symplectic structure. It does not, seemingly, preserve the Riemannian
metric nor the complex structure (the codimension of the singular set is not
 even;
see also \SS 6).

These functions $\ell,d$ lifted to $\stic$ have simple matrix-theoretic
interpretations. For $(a,b) \in \stic$, $i=1,\ldots,m$,
introduce the truncated matrices
$M_i = \pmatrix{a_1 & b_1\cr \vdots & \vdots \cr a_i & b_i\cr}$,
the first $i$ rows of $(a,b)$.
Then the $2\times 2$ matrix
$$ M_i^* M_i =  \sum_{j=1}^i \pmatrix{        |a_j|^2 & \bar{a_j} b_j \cr
                                        a_j \bar{b_j} & |b_j|^2       \cr}
$$
has the eigenvalues
$$ {1\over 2} \Bigg( \sum_{j=1}^i (|a_j|^2 + |b_j|^2) \pm
     \sqrt{   \Bigg( \sum_{j=1}^i (|a_j|^2 - |b_j|^2) \Bigg)^2
                   +\bigg|\sum_{j=1}^i a_j \bar{b_j} \bigg|^2 }\,\,\Bigg).$$

These are calculable from $\ell$ and $d$, since
$$
\eqalign{ \ell(\Phi(a,b))
&= \ell( \ldots, \phi(a_i,b_i), \ldots ) \cr
&= ( \ldots, |a_i|^2 + |b_i|^2, \ldots) \cr}
$$
and
$$
\eqalign{ d(\Phi(a,b))
&= (\ldots, \bigg|\sum_{j=1}^i \phi(a_j,b_j)\bigg|, \ldots) \cr
&= (\ldots, \sqrt{   \Bigg(\sum_{j=1}^i (|a_j|^2 - |b_j|^2)\Bigg)^2
                   + \big|\sum_{j=1}^i a_j \bar{b_j} \big|^2 }, \ldots)
\cr}
$$
So $\sum_{j=1}^i \ell_j$ is the sum of the two eigenvalues of $M_i^* M_i$,
whereas $d_i$ is the difference. (Note that $\ell_1 = d_1$ as promised;
$M_1^* M_1$'s lesser eigenvalue is $0$.)

This $(2\times 2)$-matrix $M_i^* M_i$ has the same nonzero eigenvalues as
the $i\times i$ matrix $M_i M_i^*$. The latter matrix
is more relevant in that it is the upper left $i\times i$ submatrix of
the $m\times m$ matrix $(a,b) (a,b)^*$ introduced in section (3.11).

This family of Hamiltonians -- the eigenvalues of the upper left
submatrices --
has been studied already in [Th] and is called the classical
Gel$'$fand-Cetlin system (our main reference is [GS1]).
The linear relations established above between them and $d,\ell$
are summed up in the following

\proc{Theorem}
The bending flows on ${\pol{{}}{{}}}_{\alpha}$ are the residual torus action
from the Gel$'$fand-Cetlin system on the Grassmannian $\grc$.  \endproc

The Gel$'$fand-Cetlin action on the flag manifold has always been rather
mysterious (at least to us); it is pleasant that in this case it has
 a natural geometric interpretation.

The Gel$'$fand-Cetlin functions $\{ e_{ij} \}_{j \leq i}$ (the $j$th
eigenvalue of the upper left $i\times i$ submatrix) satisfy some linear
inequalities that can be established using the minimax description of
eigenvalues [Fr, p. 149].
$$ e_{i,j} \leq e_{i-1,j+1} \leq e_{i,j+1} $$
For the polygon space functions $l,d$ most of these say $0\leq 0$;
for each $i=0,\ldots,n-1$ the nontrivial inequalities are
$$ 0 \leq -d_i + \sum_{\iota=1}^i \ell_\iota
     \leq -d_{i+1} + \sum_{\iota=1}^{i+1} \ell_\iota
     \leq d_i + \sum_{\iota=1}^i \ell_\iota
     \leq d_{i+1} + \sum_{\iota=1}^{i+1} \ell_\iota. $$
But these are transparent in our situation, as they are just
the triangle inequalities!
\setbox30=\hbox{\forcote \kern 1 truecm } \count46=\enno
\setbox31=\hbox to \wd 30{{}}
$$\copy31 \ell_{i+1} \leq d_i + d_{i+1} $$
$$\box30  d_i \leq \ell_{i+1} + d_{i+1}  $$
$$\copy31 d_{i+1} \leq \ell_{i+1} + d_i $$

(The first one, $d_i \leq \sum_{\iota=1}^i \ell_\iota$,
can be proved inductively from the others starting from $d_0=0$.)

In [GS1] it is left as an exercise to show that \cit{46} are the
{\it only} inequalities satisfied, equivalently, that every point in
the convex polytope $\Gamma_m\subset \bbr^m\times \bbr^m$ defined by them
(and $d_0=d_m=0$ and $\sum_i \ell_i = 2$) is realized by some Hermitian
matrix. We show this directly:

\proc{Theorem} The image of $\pol{m}{k\geq 2}$ under the
map $(\ell,d)$ is the whole polytope $\Gamma_m$. \endproc
\count48=\enno

Proof. We construct the
polygons directly, vertex by vertex -- really establishing that the
space $\tpol{m}{k}(\alpha,\delta)$ is nonempty (and so its quotient
by $SO(k)$ is as well). We must place
each new vertex on the intersection of two $S^{k-1}$'s, one of radius
$d_{i+1}$ from the origin, the other of radius $\ell_{i+1}$ from the
previous vertex.
The inequalities $\ell_{i+1} \leq d_i + d_{i+1}$ and
$d_{i+1} \leq \ell_{i+1} + d_i$
rule out one $S^{k-1}$ containing the other; the third inequality
$d_i \leq \ell_{i+1}+d_{i+1}$ rules out their being separated balls.
So they intersect in an $S^{k-2}$, a point or the whole $S^{k-1}$,
anywhere on which we may place the new vertex. \cqfd

\cote {\bf Remarks: }\
\decale{1)} While the map $\ell$ is equivariant with respect to the usual action
 of
$S_m$ on $\bbr^m$ , the map $d$ can only be made equivariant under the
 involution
$[i \leftrightarrow (n-i)]$ , and the polytope  $\Gamma _m$ is correspondingly
 less
symmetric than the hypersimplex $\Xi_m$.

\decale{2)} That the image of $(\ell,d)$ is the same when restricted to
planar polygons has the flavor of a more general theorem of Duistermaat [D]
on restricting moment maps to the fixed-point sets of antisymplectic
involutions.
In fact Duistermaat's theorem does not apply directly, because the subset
where $d$ is smooth (and a moment map) is noncompact; in any case we
preferred to give a polygon-theoretic proof.

\decale{3)} When $k=3$ Theorem \cit{47} guarantees that the bending torus acts
simply transitively on the fiber over an interior point of $\Gamma_m$,
making  this fiber  a torus $U(1)^{m-3}$
(or $O(1)^{m-3}$ when $k=2$). Over a prodigal
boundary point of $\Gamma_m$, the fiber is still a product of  0- or
1-spheres, but fewer of them.

\decale{4)}  Bending around other diagonals then the ones above can be done
in the same way, the moment map lifted to $\stic$ being the difference of the
 two
eigenvalues of $M^*M$ for a
corresponding submatrix $M$ of $(a,b)\in \stic$. For instance, we take
$$M=\pmatrix{a_2 & b_2\cr a_3 & b_3\cr a_4 & b_4\cr}$$
for the diagonal $\partial_{2,4}:=\rho (2)+\rho (3)+\rho (4)$.
The bending flows around two diagonals $\partial_{p,q}$ and $\partial_{p',q'}$
 commute
if and only if the pairs $\{p,q\}$ and $\{p',q'\}$ intersect or are unlinked in
$\bbr/m\bbz$.

\vskip 1 truecm \goodbreak

\noupar{Toric manifold structures on $\pol{{m}}{3}_+ (\alpha )$ for $m=4,5,6$ }

In this section, we study examples of $\pol{{}}{3}_+ (\alpha )\subset
 \pol{m}{3}$ such
that the $m-3$ diagonal functions $d_2,\dots ,d_{m-2} : \pol{{}}{3}_+ (\alpha
 )\fl{}{}
\bbr$ never vanish. The whole space $\pol{{}}{3}_+ (\alpha )$ consists of
prodigal polygons and, by \SS 5, the bending flows give an action of a big
(i.e. half-dimensional) torus  on  $\pol{{}}{3}_+ (\alpha )$.
By Delzant's theorem (see  [De], or [Gu, \SS 1]),
we can construct
from the moment polytope $\Delta _\alpha $ alone a toric manifold which is
equivariantly symplectomorphic to the space $\pol{{}}{3}_+ (\alpha )$.
This can be achieved also by [DJ,\SS 1.5], though only
up to equivariant diffeomorphism.
The latter also gives
the real part, the planar polygon space $\pol{{}}{2} (\alpha )$, as a
$2^{m-3}$-sheeted branched cover of $\Delta _\alpha $.  We sum up
below some results of these constructions without writing all the details.

Without explicit mention of the contrary, $\alpha $ is supposed to be generic.
Contrary to the previous sections, we do not require  that the perimeter
of our polygons is $2$. It was necessary to fix the perimeter in order to define
 the
map $\ell$ and the value $2$ is the natural choice to deal  with the map
$\Phi : \stic \fl{}{} \tpol{m}{k}$.
But $\fon{m}{k}(\alpha )$ makes sense for any $\alpha \in \bbr_{\geq 0}^m$
and so do the various moduli spaces $\pol{m}{k}(\alpha )$, etc. When
$\sum\alpha _i=2$, the polytope $\Delta _\alpha $ is a slice through the
Gel$'$fand-Cetlin moment polytope $\Gamma _m$ of \SS 5; for general $\alpha $ it
 is a
homothetic copy of this section.

\smallskip

\cote {\it $m = 4$: }\ The condition which guarantees that $d_2$ never vanishes
is $\alpha _1\ne \alpha _2$ or $\alpha _3\ne \alpha _4$. The space of
quadrilaterals  $\pol{4}{3}_+(\alpha )$ is then a compact toric manifold of
 dimension $2$,
therefore diffeomorphic to $\bbc P^1$. The moment map $d_2$ has image the
 interval
$\Delta _\alpha:=I_1 \cap I_2$ where
$$ I_1:= [|\alpha _1-\alpha _2|,\alpha _1+\alpha _2]
\quad \hbox{ and }\quad I_2:=
[|\alpha _4-\alpha _3|,\alpha _4+\alpha _3]$$
The space $\pol{4}{2}(\alpha )$ is $\bbr P ^1$. The quadrilateral spaces
$\pol{4}{2}(\alpha )_+$ have long since been classified (see for instance [Ha]).
One has $$ \pol{4}{2}(\alpha )_+ = \cases{ S^1 {\sqcup} S^1   &  when
\hbox{$ I_1\subset I_2$ or $I_2 \subset I_1$} \cr
S^1  &   otherwise \cr}$$
Observe also that $\alpha $ is generic if and only if the boundaries of the
 intervals
$I_1$ and $I_2$ do not meet.

By the Duistermaat-Heckman Theorem [Gu, \SS 2], the symplectic volume of
$\pol{4}{3}(\alpha )$ is proportional to the length of $\Delta _\alpha $. We
 would
then obtain the same length if we had used the other diagonal
 $|\rho (2)+\rho (3)|$. This produces a statement of elementary Euclidean
 geometry:
the variation intervals of the two diagonals of a quadrilateral with given
sides in $\bbr^3$ are the same length.

 \smallskip

\cote {\it $m = 5$: }\  Conditions for which both  $d_2$ and $d_3$ never vanish
are for instance $\alpha _1\ne \alpha _2$ and $\alpha _4\ne \alpha _5$. The
 space of
pentagons  $\pol{5}{3}_+(\alpha )$ is then a toric manifold of dimension $4$.
 The moment polytope $\Delta _\alpha  \in \bbr^2$ for $(d_2,d_3)$ is the
intersection of the rectangle $I_\alpha $
$$I_\alpha := [|\alpha _1-\alpha _2|,\alpha _1+\alpha _2]\times
[|\alpha _5-\alpha _4|,\alpha _5+\alpha _4]$$
with the non-compact rectangular region
$$\Omega _\alpha :=\{(x,y)\in (\bbr_{\geq 0})^2\mid x+y \geq \alpha _3\hbox{ and
 }
y \geq x-\alpha _3 \hbox{ and } y \leq x + \alpha _3\}$$
\line{\hfill PICTURE \hfill }
One sees that $\Delta _\alpha$ has at most 7 sides. The generic $\alpha $ are
exactly those for which the boundary of $\Omega _\alpha $  contains no corner
of $I_\alpha $ and
 $\pol{5}{3}_+(\alpha )$ is then obtained by symplectic blowings up from $\bbc
 P^2$ or
$S^2\times S^2$. The space of planar polygons $\pol{5}{2}_+(\alpha )$ is a
 closed
surface  obtained by gluing 4 copies of $\Delta _\alpha $ and its Euler
 characteristic
is given by the formula
$$\chi (\pol{5}{2}(\alpha )) = 4 - \hbox{\#} (\hbox{sides of }\Delta
_\alpha )$$
(see [DJ], Example 1.20) and is orientable if and only if
$I_\alpha \subset \omega _\alpha $. One has of course
$\chi (\pol{5}{2}_+(\alpha ))=2\chi (\pol{5}{2}(\alpha ))$ and is
$\pol{5}{2}_+(\alpha )$ is an orientable surface ($\pol{m}{k}_+(\alpha )$ is
 always
orientable). The possible cases, depending on the number of sides of $\Delta
 _\alpha
$,  are summed up in the following table.\smallskip

\def\gtvi{\vrule height 15pt depth 5 pt   width 0pt}
\def\trait{\noalign{\hrule}}

$$\vbox{ \offinterlineskip
\halign {\vrule\vrule\kern 5pt #\quad &\vrule\kern 7pt#\kern 5pt
&\vrule\kern 7pt#\kern 5pt &\vrule\kern 7pt#\kern 5pt&
\vrule\kern 7pt#\kern 5pt&\vrule\vrule#
\cr
\trait \trait %\tmVII
\gtvi \# of sides  & $\pol{{}}{3}_+ (\alpha )$ & $\pol{{}}{2} (\alpha )$
 & $\pol{{}}{2}_+ (\alpha )$  & %\tmVII
Ex. of $\alpha $     & \cr
\trait  \trait \tvi \tm
{\hfill 3\hfill }  & $\bbc P^2$ & $\bbr P^2$
 & $S^2$ & (2,1,5,1,2) & \cr \noalign{\hrule} \tvi
  &a)\kern 3pt $\bbc P^2 \# \overline{\bbc P^2}$ & Klein bottle
 & $T^2$ & (3,2,5,1,2) & \cr
\tvi {\hfill 4\hfill }  & {\hfil or \hfil\hfil } &  &  &  & \cr
 & b) {$S^2\times S^2$\hfill } & $T^2$  & $T^2 \sqcup T^2$ & (3,1,3,1,3) & \cr
 \trait
\tvi {\hfill 5\hfill }  & $(S^2\times S^2) \# \overline{\bbc P^2}$
 & $T^2\# \bbr P^2$  & $\Sigma_2$ & (2,1,3,1,2) & \cr \trait
\tvi {\hfill 6\hfill }  & $(S^2\times S^2) \# 2\overline{\bbc P^2}$
 & $T^2\# 2\bbr P^2$  & $\Sigma_3$ & (4,2,2,2,4) & \cr \trait
\tvi {\hfill 7\hfill }  & $(S^2\times S^2) \# 3\overline{\bbc P^2}$
 & $T^2\# 3\bbr P^2$  & $\Sigma_4$ & (4,2,4,2,4) & \cr \trait
}\hrule\hrule}
$$

\cote Some embeddings of
the regular pentagon $\alpha =(1,1,1,1,1)$ are not prodigal. However none
are lined and thus the moduli space $V_0 :=\pol{5}{3}(\alpha )$
is diffeomorphic for small $\varepsilon $  to $V_\varepsilon $ where
$V_\varepsilon := \pol{5}{3}(\alpha _\varepsilon )$ and
$\alpha _\varepsilon :=(1+\varepsilon,1,1,1,1+\varepsilon)$.
 The moment polytope  for $\alpha _\varepsilon $ has then 7 sides and thus
$V_0 \simeq V_\varepsilon $ is diffeomorphic to $(S^2\times S^2) \#
 3\overline{\bbc
P^2}$ (if $k=2$, $(\pol{5}{2}(\alpha )_+\simeq \Sigma _4$).
\line{\hfill PICTURE \hfill }

The pre-image in $V_\varepsilon $ of the segments $\{x=\varepsilon \}\cap \Delta
_\alpha '$ and  $\{y=\varepsilon \}\cap \Delta _\alpha '$ are 2-spheres of
symplectic volume proportional to $\varepsilon $, by the Duistermaat-Heckman
 Theorem.
Passing to the limit $V_0$, these spheres become Lagrangian, and so cannot be
complex. This shows that the action of the bending torus is not complex --
these polygon spaces are only equivariantly symplectomorphic,
not equivariantly isometric, to toric varieties.

\smallskip

\cote  \count54=\enno Any  class $r \in
\pol{5}{{k=2,3}}(\alpha )$ has a unique representative in  $\rho \in
 \tpol{5}{k}(\alpha
)$ with $\rho (5) =(-\alpha _5,0,0)$ and  $\gamma(r):=\rho (1)+\rho (2)$ in the
half-plane $\calh =\{z=0, y\geq 0\}$. This provides a map $\gamma :
 \pol{5}{3}(\alpha
)\fl{}{} \cal H$ whose image $\tilde \Delta _\alpha $  is the intersection
$R_1\cap R_2\cap \cal H$ where $R_1$ and $R_2$ are the rings
$$R_1:=\{v\in \bbr^2\mid |\alpha _1-\alpha -2| \leq |v| \leq \alpha _1+\alpha
 _2\}
\quad
R_2:=\{v\in \bbr^2\mid |\alpha _4-\alpha -3| \leq |v| \leq \alpha _4+\alpha
 _3\}$$

\picture

The idea of reconstructing $\pol{5}{{2}}(\alpha )$ by gluing copies of
$\tilde \Delta _\alpha $ goes back to the early works of W. Thurston on planar
linkages (see [Th, p.100]). The relationship with our theory is the
following: the domain $\tilde \Delta _\alpha $ is straightened up
into a PL-polytope $\Delta _\alpha$ in  $\bbr^2$ by the map $v\mapsto
(|v|,|v-(0,\alpha _5)|)$ and  $\Delta _\alpha$ is just the moment polytope for
 the
bending Hamiltonians  $\partial _1(\rho ) = |\rho (1)+ \rho (2)|$ and $\partial
_2(\rho ) = |\rho (3)+ \rho (4)|$.

\cote {\it $m = 6$: } \count55=\enno  \ The conditions $\alpha _1\ne \alpha _2$
 and
$\alpha _5\ne \alpha _6$ imply that $d_2$ and $d_4$ never vanish. However, one
 cannot
guarantee generically  $d_3\ne 0$. But we can replace
the  $d=(d_1,d_2,d_3)$ by $\delta :=(\partial  _1,\partial  _2,\partial _3)$
 where
$$\partial _1 := d_1= |\rho (1)+ \rho (2)|\ , \
\partial _2 :=  |\rho (3)+ \rho (4)|\ , \
\partial _3 := d_3= |\rho (5)+ \rho (6)|$$
and guarantee  non-vanishing of the $\delta _i$'s by the generic condition
$\alpha _{2i-1}\ne
\alpha _{2i}$. Observe that $\partial _i\pcirc \Phi :\stic \fl{}{} \bbr$
 ($i=1,2,3$) are
the functions on $\stic$ given (on $(a,b) \in \stic$) by the difference of the
eigenvalues of the $(2\times 2)$-matrices $$ \pmatrix{a_1 & b_1\cr a_2 & b_2\cr}
 \quad \quad \pmatrix{a_3 & b_3\cr a_4 & b_4\cr} \quad \quad
\pmatrix{a_5 & b_5\cr a_6 & b_6\cr}$$
The moment polytope in $\bbr^3$ is the intersection of  the rectangular
parallepiped
$$I_\alpha := [|\alpha _1-\alpha _2|,\alpha _1+\alpha _2]\times
[|\alpha _4-\alpha _3|,\alpha _4+\alpha _3]\times
[|\alpha _6-\alpha _5|,\alpha _6+\alpha _5]$$
with the region
$$\Omega  :=\{(x,y,z) \in \bbr^3 \mid
0 \leq z \leq x+y  \ , \ 0 \leq x \leq y+z \hbox{ and } 0 \leq y \leq x + z\}$$
The domain $\Omega$ can be described as the convex hull of the three half-lines
$$\{0 \leq x=y \hbox{ and } z = 0\} \ , \
\{0 \leq y=z \hbox{ and } x = 0\} \ , \
\{0 \leq z=x \hbox{ and } y = 0\}$$
or the cone $\bbr_+\cdot \Xi_3$ on the hypersimplex $\Xi_3$. The polytope
$\Delta _\alpha $ has then at most $9$ facets.
The length-system $\alpha $ is generic when the boundary of $\Omega $ do not
 contain
corners of $I_\alpha $. As 6 is even, the regular hexagon is not generic:
$\pol{6}{1}(1,\dots,1)$ contains 10 elements.

\smallskip

\cote  The bending flows $\partial $ occuring in \cit{54} and \cit{55} admit the
following generalization.
For $m=2n-1$ or $2n$, we define the {\it even-step} map
$e : \fon{m}{k}\fl{}{} \fon{n}{k}$ by $e(\rho )(i) := \rho (2i-1)+\rho (2i)$
taking $e(\rho )(n) := \rho (m)$ if $m$ is odd. We also call $e$ the induced
 maps
$\tpol{m}{k}\fl{e}{} \tpol{n}{k}$, $\pol{m}{k}_+\fl{e}{} \pol{n}{k}_+$ and
$\pol{m}{k}\fl{e}{} \pol{n}{k}$. We call $\rho \in \fon{m}{k}$ {\it even
 generic} if
$e(\rho )$ is a proper polygon. Above the space of proper polygons, the map $e$
 is a
smooth locally trivial bundle whose fiber is a product of $(k-1)$-spheres.

\picture

Define $\partial= (\partial _1,\dots , \partial _n) : \fon{m}{k}\fl{}{} \bbr^n$
 by
$\partial := \ell \pcirc e$. The  map $\partial $ gives the side lengths of
the new polygon $e(\rho )$.
It is always continuous and smooth when $e(\rho )$ is a proper
polygon.  As the map $e$ is a submersion on even-generic polygons, the critical
 values
of $\partial $ are the same as those of $\ell$, the walls of
\ecit{4}{38}. As for the map $\ell$, the
map $\partial $ can be defined on each $\pol{m}{k}(\alpha)$.
Call $\alpha \in \bbr^m$ {\it even generic} if $\pol{m}{k}(\alpha)$ only
 consists
 of even-generic polygons. For instance, $\alpha $ is even-generic if
$\alpha _{2i-1}\ne \alpha _{2i}$ for all $i$.
When $k=3$,  $\partial $
is a moment map for the corresponding bending action of $T^{n}$ defined on
even-generic polygons.

Restrict to $\pol{m}{3}(\alpha)_+$ for an even-generic
$\alpha $.
Define the right-angled polytope
$$I_\alpha := \prod_{i=1}^n [|\alpha _{2i}-\alpha _{2i-1}|,
\alpha _{2i}+\alpha _{2i-1}]$$
and consider the convex polytope $\Delta _\alpha \subset \bbr^n$
$$\Delta _\alpha := \cases{I_\alpha \cap (\bbr_+\cdot \Xi_n)  & when
$m=2n$ \cr
I_\alpha \cap (\bbr_+\cdot \Xi_n) \cap \{x_n=|\rho (m)|\}  &  when $m=2n-1$\cr
}$$

\proc{Proposition} 1) \ The image of
 $\partial : \pol{m}{k}(\alpha)_+\fl{}{} \bbr^n$
 is the whole polytope $\Delta _\alpha$.
\decale{2)} If $x \in \Delta _\alpha $ is a regular value of $\partial $, the
 even-step
map $e$ induces, for $m=3$, a symplectomorphism from the symplectic reduction
$T^n\backslash \partial ^{-1}(x)$ onto  $\pol{n}{k}_+(x)$.
\cqfd \endproc \count53=\enno

\vskip 1 truecm \goodbreak

\noupar{Remarks and open problems}   %%%%%% SECTION 7

\cote Is there an octonionic version of Section 3?
Alternately, are there $U_1(\bbh)$ bendings in
dimension 5 (like the $U_1(\bbc)$ bending flows in dimension 3 and
$U_1(R)$ flippings in dimension 2)?

\cote Observe that the inclusion $\pol{m}{k} \subset \pol{m}{k+1}$ becomes
a bijection when $k\geq m-1$ (triangles are always planar, etc.).
In what ways are these spaces $\pol{m}{m-1}$ more natural than the unstable
 ones?

\cote The $m$-polygons whose first diagonal is of a given length forms a sphere
 bundle
over a space of $(m-1)$-polygons. (For $k=3$ this is just symplectic reduction
by the first bending circle.) This gives an inductive way to construct
the space of $m$-polygons by gluing together (sphere bundles over) the spaces of
$(m-1)$-polygons; it would require identification of these
sphere bundles, which in $k=3$ might be done using the Duistermaat-Heckman
theorem (where the circle bundle is determined by its Euler class).

Alternately one might work out the fibers of the whole map $d$ of section 5.
Unfortunately in dimensions above 3 these are always singular (at, in
 particular, the planar
polygons).

\cote In [KM1] and [Wa] there are presented ``wall-crossing arguments'' for
identifying the spaces $\pol{m}{2}(\alpha)$. It would be nice to
relate these to a combination of [Du] and the paper [GS2], which presents
its own wall-crossing arguments for any symplectic reduction
by a torus.

\cote A space of great interest nowadays is the moduli space
of flat $SU(2)$ connections on a punctured Riemann sphere ---
in the language of this paper, geodesic
polygons in $S^3$ (rather than $\bbr^3$). The spaces here can be seen
as limiting versions where the radius of $S^3$ goes
to infinity. We do not know how to adapt the Gel$'$fand-MacPherson
correspondence to this case; one definite complication is that it is
no longer the symmetric group but the braid group which permutes the edges,
and that action is not complex.

\cote By averaging the Riemannian metric with respect to the bending
torus, one can deform the complex
structure on a space of prodigal polygons to that of the corresponding toric
 variety.
Is the original complex structure that of a toric variety (not just in the same
deformation class)?

\vskip .7 truecm \goodbreak

{\bf REFERENCES} \par  \vskip .4 truecm

\ref
\def\key{Au}
\def\auteur{M. Audin}
\def\titre{The topology of torus actions on symplectic manifolds.}
\def\journal{Birkh\"auser}
\def\numero{} \def\annee{(1991)}  \def\pp{}
\endref

\ref
\def\key{DJ}
\def\auteur{M. Davis \& T. Januszkiewicz}
\def\titre{Convex polytopes, Coxeter orbifolds and torus actions.}
\def\journal{Duke Math. J.}
\def\numero{62} \def\annee{(1991)}  \def\pp{417--451}
\endref

\ref
\def\key{De}
\def\auteur{T. Delzant}
\def\titre{Hamiltoniens p\'eriodiques et image convexe de l'application moment.}
\def\journal{Bull. Soc. Math. France}
\def\numero{116} \def\annee{(1988)}  \def\pp{315--339}
\endref

\ref
\def\key{Du}
\def\auteur{J.J.  Duistermaat}
\def\titre{Convexity and tightness for restrictions of Hamiltonian functions to
fixed-point sets of an antisymplectic involution.}
\def\journal{Trans. AMS}
\def\numero{275} \def\annee{(1983)}  \def\pp{417--429}
\endref

\ref
\def\key{Fr}
\def\auteur{J. Franklin}
\def\titre{Matrix theory.}
\def\journal{Prentice--Hall, Inc.}
\def\numero{} \def\annee{(1968)}  \def\pp{}
\endref

\ref
\def\key{GGMS}
\def\auteur{I. Gel$'$fand \& M. Goresky \& R. MacPherson \& V. Serganova}
\def\titre{Combinatorial geometries, convex polyedra and Schubert cells.}
\def\journal{Adv. Math.}
\def\numero{63} \def\annee{(1987)}  \def\pp{301--316}
\endref

\ref
\def\key{GM}
\def\auteur{I. Gel$'$fand  \& R. MacPherson}
\def\titre{Geometry  in Grassmannians and a generalization of the dilogarithm.}
\def\journal{Adv. Math.}
\def\numero{44} \def\annee{(1982)}  \def\pp{279--312}
\endref

\ref
\def\key{Gu}
\def\auteur{V. Guillemin}
\def\titre{Moment maps and combinatorial invariants of Hamiltonian
 $T^n$-spaces.}
\def\journal{Birkh\"auser}
\def\numero{} \def\annee{(1994)}  \def\pp{}
\endref

\ref
\def\key{GS1}
\def\auteur{V. Guillemin and S. Sternberg}
\def\titre{The Gelfand-Cetlin system and quantization of the complex flag
 manifolds.}
\def\journal{J-Funct-Anal}
\def\numero{ 52(1)} \def\annee{(1983)}  \def\pp{106-128}
\endref

\ref
\def\key{GS2}
\def\auteur{V. Guillemin and S. Sternberg}
\def\titre{Birational equivalence in the symplectic category.}
\def\journal{Inventiones Math.}
\def\numero{97} \def\annee{(1989)}  \def\pp{485--522}
\endref

\ref
\def\key{Ha}
\def\auteur{J.-C. Hausmann}
\def\titre{Sur la topologie des bras articul\'es.}
\def\journal{In ``Algebraic Topology, Poznan", Springer Lectures Notes}
\def\numero{1474} \def\annee{(1989)}  \def\pp{146--159}
\endref

\ref
\def\key{Ho}
\def\auteur{H. Hopf}
\def\titre{\"Uber di Abbildungen der dreidimensionalen Sph\"are auf die
Kugelfl\"ache} \def\journal{Math. Annalen}
\def\numero{104} \def\annee{(1931)}  \def\pp{637--665}
\endref

\ref
\def\key{KM1}
\def\auteur{M. Kapovich \& J. Millson}
\def\titre{On the moduli space of polygons in the Euclidean plane.}
\def\journal{To appear in J. of Diff. Geometry}
\def\numero{} \def\annee{}  \def\pp{}
\endref

\ref
\def\key{KM2}
\def\auteur{M. Kapovich \& J. Millson}
\def\titre{The symplectic geometry of polygons in Euclidean space.}
\def\journal{Preprint}
\def\numero{} \def\annee{(Dec. 1994)}  \def\pp{}
\endref

\ref
\def\key{Ki}
\def\auteur{F. Kirwan}
\def\titre{Cohomology of quotients in symplectic and algebraic geometry.}
\def\journal{Princeton University Press}
\def\numero{} \def\annee{(1984)}  \def\pp{}
\endref

\ref
\def\key{Mo}
\def\auteur{R. Montgomery}
\def\titre{Heisenberg and isoholonomic inequalities: ``How to use dual
pairs to ignore the rest of the universe''}
\def\journal{Symplectic geometry and mathematical physics. Souriau Volume.
Birkh\"auser } \def\numero{} \def\annee{(1991)}  \def\pp{}
\endref

\ref
\def\key{Th}
\def\auteur{A. Thimm}
\def\titre{Integrable geodesic flows on homogeneous spaces.}
\def\journal{Ergodic theory and dynamical systems}
\def\numero{1} \def\annee{(1980)}  \def\pp{495--517}
\endref

\ref
\def\key{TW}
\def\auteur{W. Thurston \& J. Weeks}
\def\titre{The mathematics of three-dimensional manifolds.}
\def\journal{Scientific American}
\def\numero{251} \def\annee{(1984)}  \def\pp{94--107}
\endref

\ref
\def\key{Wa}
\def\auteur{K. Walker}
\def\titre{Configuration spaces of linkages}
\def\journal{Undergr. Thesis, Princeton}
\def\numero{} \def\annee{(1985)}  \def\pp{}
\endref

\end